\documentclass[12pt,leqno]{article}
\usepackage{amssymb,graphics,latexsym}
\usepackage{graphicx}
\usepackage{amsmath}
\def\R{\hbox{\bf\rlap{I}{\hbox to 2 pt{}}R}}

\newcommand{\re}{{\rm Re\, }}
\newcommand{\im}{{\rm Im\, }}

\newcommand{\dia}{{\rm diag\, }}

\begin{document}
\thispagestyle{empty}
\begin{center}
\section*{Weighted Shift Matrices: Unitary Equivalence, Reducibility and Numerical Ranges}

{\bf Hwa-Long Gau$^{\mbox{a}}$,\hspace*{.3cm}\bf Ming-Cheng
Tsai$^{\mbox{b}}$,\hspace*{.3cm}\bf Han-Chun Wang$^{\mbox{a}}$}
\end{center}

\noindent ${ }^{\mbox{a}}$Department of Mathematics, National
Central
University, Chungli 32001, Taiwan\\
${ }^{\mbox{b}}$Department of Applied Mathematics, National Sun
Yat-sen University, Kaohsiung \\ \hspace*{1.5mm} 804, Taiwan

\vspace{2mm}

\noindent {\bf Abstract.}

An $n$-by-$n$ ($n\ge 3$) weighted shift matrix $A$ is one of the form
$$\left[\begin{array}{cccc}0 & a_1 & & \\ & 0 & \ddots & \\ & & \ddots & a_{n-1} \\ a_n & & & 0\end{array}\right],$$
where the $a_j$'s, called the weights of $A$, are complex numbers. Assume that all $a_j$'s are nonzero and $B$ is an $n$-by-$n$ weighted shift matrix with weights $b_1, \ldots, b_n$. We show that $B$ is unitarily equivalent to $A$ if and only if $b_1\cdots b_n=a_1\cdots a_n$ and, for some fixed $k$, $1\le k \le n$, $|b_j| = |a_{k+j}|$ ($a_{n+j}\equiv a_j$) for all $j$. Next, we show that $A$ is reducible if and only if $A$ has periodic weights, that is, for some fixed $k$, $1\le k \le \lfloor n/2\rfloor$, $n$ is divisible by $k$, and $|a_j|=|a_{k+j}|$ for all $1\le j\le n-k$. Finally, we prove that $A$ and $B$ have the same numerical range if and only if $a_1\cdots a_n=b_1\cdots b_n$ and $S_r(|a_1|^2, \ldots, |a_n|^2)=S_r(|b_1|^2, \ldots, |b_n|^2)$ for all $1\le r\le \lfloor n/2\rfloor$, where $S_r$'s are the circularly symmetric functions.

\noindent
\emph{Mathematics subject classification} (2010): 15A60.\\
\emph{Keywords and phrases}: Numerical range, reducibility,
weighted shift matrices.

{\em Email address:} hlgau@math.ncu.edu.tw (H.-L. Gau); mctsai2@gmail.com (M.-C. Tsai); 942401005@cc.ncu.edu.tw (H.-C. Wang).

${}^{\rm a}$Research supported in part by the National Science
Council of the Republic of China under project NSC
100-2115-M-008-004.

\newpage

\centerline {\bf\large 1. Introduction}

\vspace{.5cm}

An $n$-by-$n$ ($n\ge 3$) weighted shift matrix $A$ is one of the
form
$$\left[\begin{array}{cccc}0 & a_1 & & \\ & 0 & \ddots & \\ & & \ddots & a_{n-1} \\ a_n & & & 0\end{array}\right],$$
where the $a_j$'s, called the weights of $A$, are complex numbers.
In this paper, we study some unitary-equivalence properties of
such matrices. Previous works in this respect are the necessary
and sufficient conditions for the boundary of numerical range of
$A$ to have a line segment \cite[Theorem 1]{5}. Here we consider
other properties of these matrices such as their reducibility and
their numerical range.

\vspace{.5cm}

In Section 2, we give necessary and sufficient conditions for two
$n$-by-$n$ weighted shift matrices $A$ and $B$ with weights $a_1,
\ldots , a_n$ and $b_1, \ldots , b_n$, respectively, to be
unitarily equivalent. More specifically, it is shown that if $A$ has at most two zero weights, then
$B$ is unitarily equivalent to $A$ if and only if, for some fixed $k$, $1\leq k\leq n$,
$|b_j| = |a_{k+j}|$ $(a_{n+j} \equiv a_j)$ for all $1\leq j\leq
n$, and $a_1\cdots a_n = b_1\cdots b_n$. In Section 3 below, we
solve the problem when a weighted shift matrix is reducible, that
is, when it is unitarily equivalent to the direct sum of two other
matrices. We obtain a complete characterization of reducibility in
terms of the weights. It roughly says that a weighted shift matrix
is reducible when it has at least two zero weights or its weights are periodic. We
take up the numerical ranges of weighted shift matrices in Section
4. We have known that the numerical range of an $n$-by-$n$ matrix
$A$ is completely determined by its Kippenhahn polynomial
$p_A(x,y,z) = \det (x\re A+y\im A+zI_n)$, where $\re A =
(A+A^*)/2$ and $\im A = (A-A^*)/(2i)$ are the real and the
imaginary part of $A$, respectively, and $I_n$ denotes the
$n$-by-$n$ identity matrix (cf. \cite[Theorem 10]{7}). We give an
explicit expansion of $p_A(x,y,z)$ in terms of the weights of an
$n$-by-$n$ weighted shift matrix $A$. Finally, let $A$ and $B$ be
$n$-by-$n$ weighted shift matrices with weights $a_1, \ldots ,
a_n$ and $b_1, \ldots , b_n$, respectively, we give necessary and sufficient
conditions for $A$ and $B$ to have the same numerical range. More
specifically, it is shown that the following statements are
equivalent: (a) $W(A) = W(B)$; (b) $p_A(x,y,z) = p_B(x,y,z)$; (c)
$a_1\cdots a_n = b_1\cdots b_n$ and $S_{r}(|a_{1}|^2,\ldots
,|a_{n}|^2) = S_{r}(|b_{1}|^2,\ldots ,|b_{n}|^2)$ for all $1\leq
r\leq \lfloor n/2\rfloor$, where $S_r$'s are the circularly
symmetric functions (see \cite[P. 496]{6}) for more details).

\vspace{.5cm}

For any
nonzero complex number $z = x + iy$ ($x$ and $y$ real), $\arg z$
is the angle $\theta$, $0\le\theta< 2\pi$, from the positive
$x$-axis to the vector $(x, y)$. For an $n$-by-$n$ matrix $A$, let $A^*$ denote its adjoint and
$\sigma (A)$ its spectrum. Throughout this paper,
if $a_1, \ldots , a_n$ are the weights of an $n$-by-$n$ weighted shift
matrix, we always assume that
$a_{n+j}\equiv a_j$ and $a_{j-n}\equiv a_j$ for all $1\leq j\leq n$.

\vspace{1cm}

\centerline {\bf\large 2. Unitary equivalence}

\vspace{.5cm}

In \cite{5}, the authors gave sufficient conditions for unitary equivalence of two $n$-by-$n$ weighted shift matrices (cf. \cite[Lemma 2]{5}).
For the convenience of the reader we repeat this result without proofs, thus making our exposition self-contained.

\vspace{.5cm}

{\bf Lemma 2.1. \cite{5}}
\emph{Let $A$ and $B$ be $n$-by-$n$ weighted shift matrices with weights
$a_1, \ldots , a_n$ and $b_1, \ldots , b_n$}, \emph{respectively}.

(a) \emph{If}, \emph{for some fixed $k$}, $1\le k \le n$, $b_j = a_{k+j}$
($a_{n+j}\equiv a_j$) \emph{for all $j$}, \emph{then $A$ is unitarily
equivalent to $B$}.

(b) \emph{If $|a_j| = |b_j|$ for all $j$}, \emph{then $A$ is unitarily equivalent to $e^{i\psi_k}B$}, \emph{where $\psi_k = (2k\pi +\sum_{j=1}^n(\arg a_j-\arg
b_j))/n$ for $0 \le k < n$}.

\vspace{.5cm}

Let $A$ and $B$ be $n$-by-$n$
weighted shift matrices with weights $a_1, \ldots , a_n$ and $b_1, \ldots , b_n$, respectively.
It is natural to ask whether the converse of
Lemma 2.1 (a) is true. In this section, we give the
affirmative answer and show that if $A$ has at most two
zero weights, then $A$ is unitarily equivalent to $B$ if and only if,
for some fixed $k$, $1\le k \le n$, $|b_j| = |a_{k+j}|$
($a_{n+j}\equiv a_j$) for all $j$, and $a_1\cdots a_n=b_1\cdots b_n$.
We first give necessary conditions for unitary equivalence of two $n$-by-$n$ weighted shift matrices in the next proposition.

\vspace{.5cm}

{\bf Proposition 2.2.} \emph{Let $A$ and $B$ be $n$-by-$n$
weighted shift matrices with weights $a_1, \ldots , a_n$ and $b_1, \ldots , b_n$},
\emph{respectively}. \emph{If $A$ is unitarily equivalent to $B$},
\emph{then the following statements hold}.

(a) $\{|a_1|, |a_2|, \ldots, |a_n|\}=\{|b_1|, |b_2|, \ldots,
|b_n|\}$ (\emph{counting multiplicities}).

(b) $a_1\cdots a_n=b_1\cdots b_n$.

(c) \emph{If $a_l\cdots a_m\neq 0$ for some $l$ and $m$}, $1\le l\le m\le n$,
\emph{then there exists a fixed $k$}, $1\le k \le n$, \emph{such
that $|a_j|=|b_{k+j}|$} ($b_{n+j}\equiv b_j$) \emph{for} $j=l-1,
l, \ldots, m+1$.

\vspace{.5cm}

{\em Proof.} (a) A simple computation shows that $AA^*=$ diag
$(|a_1|^2, \ldots, |a_n|^2)$ and $BB^*=$ diag $(|b_1|^2, \ldots,
|b_n|^2)$. Thus the singular values of $A$ (resp., $B$) are
$|a_1|, \ldots , |a_n|$ (resp., $|b_1|, \ldots , |b_n|$ ). Since
$A$ is unitarily equivalent to $B$, $A$ and $B$ have the same
singular values, hence $\{|a_1|, |a_2|, \ldots, |a_n|\}=\{|b_1|,
|b_2|, \ldots, |b_n|\}$ (counting multiplicities) as desired.

\vspace{.5cm}

(b) An easy computation shows that $\det A = (-1)^{n+1}a_1\cdots
a_n$ and $\det B = (-1)^{n+1}b_1\cdots b_n$. Since $A$ is
unitarily equivalent to $B$, hence $\det A = \det B$ or $b_1\cdots
b_n = a_1\cdots a_n$ as asserted.

\vspace{.5cm}

(c) From Lemma 2.1 (a), we may assume that $l=1$.
By assumption, there exists an $n$-by-$n$ unitary matrix
$U=[u_{i j}]_{i,j=1}^{n}$ such that $AU = UB$. It follows that
$(A^jA^{j*})U = U(B^jB^{j*})$ for all $j$, $1\leq j\leq n$. A
direct computation shows that $A^jA^{j*} = $ diag
$(\alpha_1^{(j)}, \ldots , \alpha_n^{(j)})$ and $B^jB^{j*} = $
diag $(\beta_1^{(j)}, \ldots , \beta_n^{(j)})$, where
$\alpha_t^{(j)} = |a_ta_{t+1}\cdots a_{t+j-1}|^2$ and $\beta_t^{(j)} = |b_tb_{t+1}\cdots b_{t+j-1}|^2$
for all $1\leq t,j\leq n$. Since $U$ is
unitary, then the first row of $U$ must have a nonzero entry, that
is, $u_{1 k}\neq 0$ for some $k$, $1\leq k\leq n$. Now, we
consider the $(1,k)$-entry of $(A^jA^{j*})U$ and $U(B^jB^{j*})$,
respectively. Note that the $(1,k)$-entry of $(A^jA^{j*})U$
(resp., $U(B^jB^{j*})$) is $\alpha_1^{(j)}u_{1 k}$ (resp.,
$u_{1 k}\beta_k^{(j)}$) for $j = 1,2,\cdots ,n$. Since
$(A^jA^{j*})U = U(B^jB^{j*})$, we obtain $$|a_1a_2\cdots
a_j|^2u_{1 k} = \alpha_1^{(j)}u_{1 k} = u_{1 k}\beta_k^{(j)} =
u_{1 k}|b_kb_{k+1}\cdots b_{k+j-1}|^2 \leqno{(1)}$$ for all $j$.
For $j = 1$, since $u_{1 k} \neq 0$, we deduce that $|a_1|
= |b_k|$ from Equation (1). For $j = 2$, since $u_{1 k} \neq 0$ and $|b_k|=|a_1|\neq 0$, by Equation (1), we infer
that $|a_2| = |b_{k+1}|$. Repeating this argument gives us $|a_t|
= |b_{(k-1)+t}|$ for $t = 1, 2, \ldots ,m+1$. On the other hand,
since $A^*A=$ diag $(|a_n|^2, |a_1|^2 \ldots, |a_{n-1}|^2)$,
$B^*B=$ diag $(|b_n|^2, |b_1|^2, \ldots, |b_{n-1}|^2)$ and
$(A^*A)U = U(B^*B)$, then $|a_n|^2u_{1 k} = u_{1 k}|b_{k-1}|^2$
or $|a_n| = |b_{k-1}|$. This completes the proof.
\hfill$\blacksquare$

\vspace{.5cm}

The following theorem is our main result in this section.

\vspace{.5cm}

{\bf Theorem 2.3.} {Let $A$ and $B$ be $n$-by-$n$ weighted shift
matrices with weights $a_1, \ldots , a_n$ and $b_1, \ldots , b_n$}, \emph{respectively}.
\emph{Suppose that $A$ has at most two zero weights}. \emph{Then
$A$ is unitarily equivalent to $B$ if and only if $a_1\cdots a_n =
b_1\cdots b_n$} \emph{and}, \emph{for some fixed $k$}, $1\le k \le
n$, $|a_j| = |b_{k+j}|$ ($b_{n+j}\equiv b_j$) \emph{for all} $j$,
$1\leq j\leq n$.

\vspace{.5cm}

{\em Proof.} The sufficiency is a consequence of  Lemma
2.1 (a) and (b), we only need prove the necessity. Clearly, we have $a_1\cdots a_n=b_1\cdots b_n$ from Proposition 2.2 (b).
If all $a_j$'s are nonzero, that is, $a_1\cdots a_n\neq 0$, our assertion follows from
Proposition 2.2 (c). If $A$ has exactly one zero weight, by
Lemma 2.1 (a), we may assume $a_n = 0$ and $a_1\cdots
a_{n-1} \neq 0$. By Proposition 2.2 (c), we have $|a_j|=|b_{k+j}|$ for all $1\leq j \leq n$, for some fixed $k$,
$1\leq k \leq n$ and we are done.

\vspace{.5cm}

Now, if $A$ has exactly two zero weights, by Lemma 2.1 (a), we may
assume that $a_1\cdots a_{m-1}\neq 0$, $a_{m+1}\cdots a_{n-1} \neq
0$ and $a_m = a_n = 0$ for some $m$, $1<m<n$. Since $a_1\cdots a_{m-1}\neq 0$, Proposition 2.2
(c) implies that, for some fixed $k$, $1\leq k\leq n$, $|a_j| =
|b_{k+j}|$ for $j = 1, 2, \ldots ,m-1$, $a_n = b_k = 0$ and $a_m =
b_{k+m} = 0$. For convenience, we let $b_j'=b_{k+j}$ for all $1\le j\le n$,
and $B'$ be the $n$-by-$n$ weighted shift matrix with weights $b_1', \ldots, b_n'$.
Lemma 2.1 (a) yields that $B'$ is unitarily equivalent to $A$. Moreover, we have
$|b_j'|=|b_{k+j}|=|a_j|$ for $j=n, 1, \ldots, m$. Note that $b_m'=a_m=0=a_n=b_n'$.
Now, we need only check that $|b_j'|=|a_j|$ for $j=m+1, \ldots, n-1$.
Indeed, since $a_{m+1}\cdots a_{n-1} \neq 0$,
by Proposition 2.2 (c), there exists a fixed $s$, $1\le s\le n$, such that $|a_j|=|b_{s+j}'|$ for $j=m, m+1, \ldots, n$.
In particular, $b_s'=b_{s+n}'=a_n=0$ and $b_{s+m}'=a_m=0$.
Note that $A$ has exactly two zero weights $a_m$ and $a_n$,
by Proposition 2.2 (a), there are also exactly two weights of $B$ being zero. It forces that either $s=n$ or $s=m$ and $s+m=n$.
If $s=n$ then $|a_j|=|b_{j}'|$ for all $m\le j\le n$ and we are done. For the latter case, it implies that $s=m=n/2$. Consequently, we have actually proved that $a_{n/2}=a_n=b_n'=b_{n/2}'=0$ and, for each $j=1, \ldots, n/2$,  $|a_{(n/2)+j}|=|b_{s+(n/2)+j}'|=|b_{n+j}'|=|b_j'|=|a_{j}|$. On the other hand, since $b_{j}'\neq 0$ for all  $j=(n/2)+1, \ldots, n-1$, Proposition 2.2 (c) yields that there exists a fixed $t$, $1\le t\le n$, such that $|b_{j}'|=|a_{t+j}|$ for $j=n/2, \ldots, n$. It implies that $a_{t+(n/2)}=b_{n/2}'=0=b_n'=a_{t+n}=a_t$. But the zero weights of $A$ are exactly $a_{n/2}$ and $a_n$, we infer that either $t=n$ or $t=n/2$. If $t=n$ then  $|b_{j}'|=|a_{n+j}|=|a_j|$ for $j=n/2, \ldots, n$, as asserted. For the latter case, $t=n/2$ implies that $|b_{j}'|=|a_{(n/2)+j}|=|a_j|$ for $j=n/2, \ldots, n$. This completes the proof. \hfill$\blacksquare$

\vspace{.4cm}

We remark that if $A$ has more than  two zero weights, then the
necessity of Theorem 2.2 is not true in general. For example, let
$A$ be the $6$-by-$6$ weighted shift matrix with weights $1, 0, 2,
0, 3, 0$. Then $A = A_1\oplus A_2\oplus A_3$, where $A_j =
\left[\begin{array}{cccc}0 & j \\ 0 & 0\end{array}\right]$ for $j
= 1, 2, 3$. Let $B = A_1\oplus A_3\oplus A_2$, then $B$ is the
$6$-by-$6$ weighted shift matrix with weights $1, 0, 3, 0, 2, 0$.
It is obvious that $A$ is unitarily equivalent to $B$, but there
is no any $k$, $1\leq k\leq 6$, such that $|a_j| = |b_{k+j}|$ for
all $j = 1,\ldots , 6$.

\vspace{.4cm}

We now consider the case of weighted shift matrices with at least three zero weights. Let $A$ be an $n$-by-$n$ weighted shift matrices with weights $a_1, \ldots, a_n$. If $A$ has at least three zero weights, by Lemma 2.1 (a) and (b), we may assume that $a_n=0$ and $a_j\ge 0$ for all $j$. In this case, $A=A_1\oplus A_2\oplus\cdots\oplus A_m$, where $A_i$ is a $k_i$-by-$k_i$ weighted shift matrices with exactly one zero weight for all $1\le i\le m$, $\sum_{i=1}^mk_i=n$ and $3\le m\le n$. Note that if $k_i=1$ then $A_i=[0]$, moreover, $A_i$ is irreducible for all $1\le i\le m$. Now, let $\tau:\{1, \ldots, m\}\rightarrow\{1, \ldots, m\}$ be a permutation and $B=A_{\tau(1)}\oplus\cdots\oplus A_{\tau(m)}$, we known that $B$ is unitarily equivalent to $A$. The next theorem shows that its converse is also true, more precisely, if $B$ be an $n$-by-$n$ weighted shift matrices with nonnegative weights and $B$ is unitarily equivalent to $A$, then $B=A_{\tau(1)}\oplus\cdots\oplus A_{\tau(m)}$ for some permutation $\tau:\{1, \ldots, m\}\rightarrow\{1, \ldots, m\}$.

\vspace{.4cm}

Here, for any $n$-by-$n$ matrix $A=[a_{ij}]_{i,j=1}^n$, let $|A|$ denote the $n$-by-$n$ nonnegative matrix $[|a_{ij}|]_{i,j=1}^n$.

\vspace{.4cm}

{\bf Theorem 2.4.} \emph{Let} $A=A_1\oplus A_2\oplus\cdots\oplus A_m$, \emph{where $A_i$ is a $k_i$-by-$k_i$ weighted shift matrices with weights} $a_1^{(i)}, \ldots,  a_{k_i-1}^{(i)}, 0$ \emph{and} $a_j^{(i)}\neq 0$ \emph{for all} $1\le j\le k_i-1$ \emph{and} $1\le i\le m$, $\sum_{i=1}^mk_i=n$ \emph{and} $3\le m\le n$. \emph{Let $B$ be an $n$-by-$n$ weighted shift matrix. Then $B$ is unitarily equivalent to $A$ if and only if $|B|=|A_{\tau(1)}|\oplus |A_{\tau(2)}|\oplus\cdots\oplus |A_{\tau(m)}|$ for some permutation $\tau:\{1, \ldots, m\}\rightarrow\{1, \ldots, m\}$}.

\vspace{.4cm}

For the proof of Theorem 2.4, we need the following lemma.

\vspace{.4cm}

{\bf Lemma 2.5.}  \emph{Let} $A$ (\emph{resp}., $B$) \emph{be an $m$-by-$m$} (\emph{resp}., $n$-\emph{by}-$n$) \emph{weighted shift matrix with weights} $a_1, \ldots,  a_{m-1}, 0$ (\emph{resp}., $b_1, \ldots,  b_{n-1}, 0$), \emph{where $a_i, b_j>0$ for all} $i, j$. \emph{If $U$ is an $m$-by-$n$ matrix so that $A^kU=UB^k$ and ${A^*}^kU=U{B^*}^k$ for all $k\ge 0$}, \emph{then the following statements hold}.

(a) \emph{If $m\neq n$ then $U=0$}.

(b) \emph{If $m=n$ then $U$ is a diagonal matrix}.

(c) \emph{If $m=n$ and $A\neq B$}, \emph{then} $U=0$.

(d) \emph{If $m=n$ and $A=B$}, \emph{then $U=\alpha I_n$ for some $\alpha\in \mathbb{C}$}.

\vspace{.4cm}

{\em Proof}. By assumption, we have $(A^k{A^*}^k)U=U(B^k{B^*}^k)$ and $({A^*}^kA^k)U=U({B^*}^kB^k)$ for all $k\ge 0$. Let $U=[u_{ij}]$ and $r=\min\{m, n\}$. For $1\le i<j\le r$, we consider the $(i,j)$-entry of  $({A^*}^iA^i)U$ and $U({B^*}^iB^i)$, respectively. By a direct computation, we have $${A^*}^iA^i=\dia(\underbrace{0, \ldots, 0}_{i \mbox{ \scriptsize copies}}, |a_1\cdots a_i|^2, \ldots, |a_{m-i}\cdots a_{m-1}|^2)$$ and $${B^*}^iB^i=\dia(\underbrace{0, \ldots, 0}_{i \mbox{ \scriptsize copies}}, |b_1\cdots b_i|^2, \ldots, |b_{n-i}\cdots b_{n-1}|^2).$$ Thus the $(i,j)$-entry of  $({A^*}^iA^i)U$ and $U({B^*}^iB^i)$ are 0 and $u_{ij}|b_{j-i}\cdots b_{j-1}|^2$, respectively. Since $({A^*}^iA^i)U=U({B^*}^iB^i)$ and $b_{j-i}\cdots b_{j-1}\neq 0$, we deduce that $u_{ij}=0$ for all $1\le i<j\le r$. On the other hand, for $1\le j<i\le r$, the $(i,j)$-entry of  $({A^*}^jA^j)U$ and $U({B^*}^jB^j)$ are $u_{ij}|a_{i-j}\cdots a_{i-1}|^2$ and 0, respectively. Since $({A^*}^jA^j)U=U({B^*}^jB^j)$ and $a_{i-j}\cdots a_{i-1}\neq 0$, we deduce that $u_{ij}=0$ for all $1\le j<i\le r$. Therefore, if $m=n$ then $U$ is a diagonal matrix. This completes the proof of (b).

\vspace{.4cm}

For the other cases, if $r=n<m$, then $B^n=0$ and $U({B^*}^nB^n)=0$. But $${A^*}^nA^n=\dia(\underbrace{0, \ldots, 0}_{n \mbox{ \scriptsize copies}}, |a_1\cdots a_n|^2, \ldots, |a_{m-n}\cdots a_{m-1}|^2),$$ since $({A^*}^nA^n)U=U({B^*}^nB^n)=0$ and all $a_j$'s are nonzero, we obtain that $u_{ij}=0$ for all $n<i\le m$ and $1\le j\le n$. We now check that $u_{ii}=0$ for all $1\le i\le n$. Indeed, for $1\le i\le n$, by a direct computation, we have $$A^i{A^*}^i=\dia(|a_1\cdots a_i|^2, \ldots, |a_{m-i}\cdots a_{m-1}|^2, \underbrace{0, \ldots, 0}_{i \mbox{ \scriptsize copies}})$$ and $$B^i{B^*}^i=\dia(|b_1\cdots b_i|^2, \ldots, |b_{n-i}\cdots b_{n-1}|^2, \underbrace{0, \ldots, 0}_{i \mbox{ \scriptsize copies}}).$$ Then the $(n-i+1, n-i+1)$-entry of $(A^i{A^*}^i)U$ and $U(B^i{B^*}^i)$ are $|a_{n-i+1}\cdots a_{n}|^2u_{n-i+1,n-i+1}$ and 0, respectively. Since $(A^i{A^*}^i)U=U(B^i{B^*}^iB^i)$ and $a_{n-i+1}\cdots a_{n}\neq 0$, it forces that $u_{ii}=0$ for all $1\le i\le n$, hence we conclude that $U=0$. Similarly, if $r=m<n$, then $A^m=0$ and $({A^*}^mA^m)U=0$. But $${B^*}^mB^m=\dia(\underbrace{0, \ldots, 0}_{m \mbox{ \scriptsize copies}}, |b_1\cdots b_m|^2, \ldots, |b_{n-m}\cdots b_{n-1}|^2),$$ since $({A^*}^mA^m)U=U({B^*}^mB^m)=0$ and all $b_j$'s are nonzero, we obtain that $u_{ij}=0$ for all $1\le i\le m$ and $m< j\le n$. Moreover, for $1\le i\le m$,  the $(m-i+1, m-i+1)$-entry of $(A^i{A^*}^i)U$ and $U(B^i{B^*}^i)$ are 0 and $|b_{m-i+1}\cdots b_{m}|^2u_{m-i+1,m-i+1}$, respectively. Since $(A^i{A^*}^i)U=U(B^i{B^*}^iB^i)$ and $b_{m-i+1}\cdots b_{m}\neq 0$, it implies that $u_{ii}=0$ for all $1\le i\le m$, hence we also obtain that $U=0$. This completes the proof of (a).

\vspace{.4cm}

For the proof of (c) and (d), we may assume that $m=n$ and $U=\dia(u_{11}, \ldots, u_{nn})$ from  (b). For $1\le i\le n-1$, the $(i,i+1)$-entry of $AU$ and $UB$ are $a_iu_{i+1,i+1}$ and $b_iu_{ii}$, respectively. Then $AU=UB$ yields that $$u_{i+1,i+1}=\frac{b_i}{a_i}u_{ii}\leqno{(2)}$$ for all $1\le i\le n-1$. Similarly,  the $(i+1,i)$-entry of $A^*U$ and $UB^*$ are $a_iu_{ii}$ and $b_iu_{i+1,i+1}$, respectively. Then $A^*U=UB^*$ yields that $$u_{i+1,i+1}=\frac{a_i}{b_i}u_{ii}\leqno{(3)}$$ for all $1\le i\le n-1$. Combining (2) with (3) yields $u_{ii}=(a_i^2/b_i^2)u_{ii}$ for all $1\le i\le n-1$. Now, if $A\neq B$ then $a_j\neq b_j$ for some $j$. Note that $a_i, b_i>0$ for all $i$, thus $a_j^2/b_j^2\neq 1$ and $u_{jj}=0$. Consequently, we have $u_{ii}=0$ for all  $1\le i\le n$ by Equation (2). This completes the proof of (c). On the other hand, if $A=B$ then $a_i=b_i$ for all $i$, it follows that $u_{11}=u_{22}=\cdots=u_{nn}$ from Equation (2). Hence $U=u_{11}I_n$, completing the proof. \hfill$\blacksquare$

\vspace{.4cm}

We are now ready to prove Theorem 2.4.

\vspace{.4cm}

{\em Proof of Theorem} 2.4. The sufficiency follows from Lemma 2.1 (b) and a suitable permutation matrix, we only need prove the necessity. Assume that $B$ is unitarily equivalent to $A$. Since $A$ has zero weights, by Lemma 2.1 (b), we may assume that all weights of $A$ are nonnegative. After a permutation of $A_i$'s, we may assume that $1\le k_1\le k_2\le\cdots\le k_m\le n$. Note that $\sigma(A)=\{0\}$ and $\dim \ker A_i=1$ for all $1\le i\le m$, the Jordan canonical form of $A$ is $J_{k_1}\oplus J_{k_2}\oplus\cdots\oplus J_{k_m}$, where $J_{k_i}$ is the Jordan block of size $k_i$ for all $i$. This means that the number and sizes of direct summands of $A$ are completely determined by its Jordan canonical form. Since $A$ and $B$ have the same Jordan canonical form, after a permutation of direct summands of $B$, we can write $B=B_1\oplus B_2\oplus\cdots\oplus B_m$, where $B_i$ is a $k_i$-by-$k_i$ weighted shift matrices with weights $b_1^{(i)}, \ldots,  b_{k_i-1}^{(i)}, 0$ and $b_j^{(i)}\neq 0$ for all $1\le j\le k_i-1$ and $1\le i\le m$. Moreover, by Lemma 2.1 (b), we also assume $b_j^{(i)}>0$ for all $1\le j\le k_i-1$ and $1\le i\le m$.

\vspace{.4cm}

We need to check that $\{A_1, \ldots, A_m\}=\{B_1, \ldots, B_m\}$ (counting multiplicities). Indeed, let $U=[U_{ij}]_{i,j=1}^m$ be a $n$-by-$n$ unitary block matrix so that $AU=UB$, where $U_{ij}$ is a $k_i$-by-$k_j$ matrix for all $1\le i, j\le m$. For each $A_i$, $1\le i\le m$, if $A_i\neq B_j$ for all $j$, $1\le j\le m$, Lemma 2.5 (a) and (c) yield that $U_{ij}=0$ for all $1\le j\le m$. This contradicts to the fact that $U$ is unitary. Thus we deduce that $\{A_1, \ldots, A_m\}\subseteq\{B_1, \ldots, B_m\}$. Similarly, if $B_j\neq A_i$ for all $1\le i\le m$, then $U_{ij}=0$ for all $1\le i\le m$, a contradiction. Hence we conclude that $\{A_1, \ldots, A_m\}=\{B_1, \ldots, B_m\}$.

\vspace{.4cm}

Next, for each $A_i$, $1\le i\le m$, we need to count the multiplicities of $A_i$ in $A$ and $B$, respectively. If $1=k_1=\cdots=k_t$ and $k_{t+1}>1$ for some $t$, that is, $A_1=\cdots=A_t=[0]$ and $A_i\neq[0]$ for all $i>t$,  since the sizes of $A_i$ and $B_i$ are the same for all $1\le i\le m$, thus $B_1=\cdots=B_t=[0]$ and $B_i\neq[0]$ for all $i>t$. This means that $\ker A\cap\ker A^*=\ker B\cap\ker B^*=\mathbb{C}^t\oplus\{0\}\subseteq\mathbb{C}^n$. Since $AU=UB$ and $A^*U=UB^*$, then $U(\ker B\cap\ker B^*)=\ker A\cap\ker A^*$, it implies that $U=[U_{ij}]_{i,j=1}^t\oplus[U_{ij}]_{i,j=t+1}^m$, $[U_{ij}]_{i,j=1}^t$ and $[U_{ij}]_{i,j=t+1}^m$ are unitary, $(A_1\oplus\cdots\oplus A_t)[U_{ij}]_{i,j=1}^t=[U_{ij}]_{i,j=1}^t(B_1\oplus\cdots\oplus B_t)$ and $(A_{t+1}\oplus\cdots\oplus A_m)[U_{ij}]_{i,j=t+1}^m=[U_{ij}]_{i,j=t+1}^m(B_{t+1}\oplus\cdots\oplus B_m)$. Therefore, we may assume that $2\le k_1\le k_2\le \cdots\le k_m$, that is, both $A$ and $B$ have no zero direct summand.

\vspace{.4cm}

Next, after permutations of $A_i$'s and $B_i$'s, respectively, we may assume that $A_1=A_2=\cdots=A_s$, $A_i\neq A_1$ for all $i>s$, $B_1=B_2=\cdots=B_r=A_1$ and $B_j\neq A_1$ for all $j>r$. We want to show $s=r$. Indeed, Lemma 2.5 (a), (c) and (d) yield $$U=\left[\begin{array}{ccc|c} \alpha_{11}I_{k_1} & \cdots &  \alpha_{1r}I_{k_1} & \\ \vdots & & \vdots & 0\\ \alpha_{s1}I_{k_1} & \cdots &  \alpha_{sr}I_{k_1} & \\ \hline &  0 & & U'\end{array}\right],$$ where $\alpha_{ij}$'s are complex numbers and $U'$ is a ($n-sk_1$)-by-($n-rk_1$) matrix. Since $U$ is unitary, the column vectors and the row vectors of $U$ are orthonormal, respectively, it forces that $s=r$. Therefore, we conclude that $A_1=\cdots=A_s=B_1=\cdots=B_s$, $A_i, B_i\neq A_1$ for all $i>s$, $U'$ is unitary and  $(A_{s+1}\oplus\cdots\oplus A_m)U'=U'(B_{s+1}\oplus\cdots\oplus B_m)$. Repeating this argument gives us $\{A_1, \ldots, A_m\}=\{B_1, \ldots, B_m\}$ (counting multiplicities). This completes the proof. \hfill$\blacksquare$

\vspace{1cm}

\centerline {\bf\large 3. Reducibility}

\vspace{.5cm}

A matrix is \emph{reducible} if it is unitarily equivalent to the
direct sum of two other matrices. Let $A$ be an $n$-by-$n$
weighted shift matrix with weights $a_1,\ldots , a_n$. In
\cite{5}, the author shown that $A$ is reducible if and only if
one of the following cases holds: (a) $a_i = a_j = 0$ for some
$1\leq i<j\leq n$, (b) $n$ is odd, $|a_j| = |a_1| \neq 0$ for all
$1\leq j\leq n$, (c) $n$ is even, $|a_j| = |a_{j+(n/2)}| \neq 0$
for all $1\leq j\leq n/2$. Unfortunately, this result is
erroneous. For example, let $A$ be the $6$-by-$6$ weighted shift
matrix with weights $1, 2, 1, 2, 1, 2$. It is obvious that $A$
does not satisfy the conditions (a), (b) and (c). Let $U$ be the
$6$-by-$6$ unitary matrix
$$\frac{1}{\sqrt{3}}\left[\begin{array}{cccccc}1 & 0 & 1 & 0 & 1 & 0 \\ 0 & 1 & 0 & \omega & 0 & \omega^2 \\ 1 & 0 & \omega^2 & 0 & \omega^4 & 0 \\  0 & 1 & 0 & \omega^3 & 0 & \omega^6 \\ 1 & 0 & \omega^4 & 0 & \omega^8 & 0 \\  0 & 1 & 0 & \omega^5 & 0 & \omega^{10} \end{array}\right],$$
where $\omega = e^{i\pi/3}$. A direct computation shows that
$$U^*AU = \left[\begin{array}{cccc}0 & 1 \\ 2 & 0\end{array}\right]\oplus \left[\begin{array}{cccc}0 & \omega \\ 2\omega & 0\end{array}\right]\oplus \left[\begin{array}{cccc}0 & \omega^2 \\ 2\omega^2 &
0\end{array}\right].$$ Hence $A$ is reducible.

\vspace{.5cm}

In this section, we give a criterion for a weighted shift matrix
$A$ to be reducible. We will show that $A$ is reducible if and
only if one of the following cases holds: (a) $a_i = a_j = 0$ for
some $1\leq i<j\leq n$, (b) $A$ has periodic weights.

\vspace{.5cm}

We first consider the weighted shift matrices with nonzero
weights.

\vspace{.5cm}

{\bf Theorem 3.1.} \emph{Let $A$ be an $n$-by-$n$ weighted shift
matrix with nonzero weights $a_1,\ldots , a_n$. Then the following
statements are equivalent:}

(a) \emph{$A$ is reducible.}

(b) \emph{$A$ has periodic weights, that is, for some fixed $k$,
$1\le k \le \lfloor n/2\rfloor$, $n$ is divisible by $k$, and $|a_j|=|a_{k+j}|$
for all $1\le j\le n-k$}.

(c) \emph{there exists a fixed $k$, $1\leq k\leq \lfloor n/2\rfloor$, such that
$n$ is divisible by $k$ and  $A$ is unitarily equivalent to
$e^{i\theta}\left(B\oplus (\omega B)\oplus \cdots \oplus
(\omega^{(n/k)-1}B)\right)$, where $B$ is the $k$-by-$k$ weighted
shift matrix with weights $|a_1|, \ldots , |a_k|$, $\omega
=e^{2\pi i/n}$ and $\theta = \left(\sum_{i=1}^{n}\arg
a_j\right)/n$}.

\vspace{.5cm}

{\em Proof.} (a)$\Rightarrow$(b). Since $A$ is reducible, there
exists an $n$-by-$n$ orthogonal projection
$P=[p_{i j}]_{i,j=1}^{n}$ such that $PA = AP$. It follows that
$PA^* = A^*P$, consequently, $P(A^kA^{k*}) = (A^kA^{k*})P$ for all $k$.

\vspace{.5cm}

On the contrary, suppose that the weights of $A$ are not periodic.
We claim that for any $1\leq i<j\leq n$, there exists $k_0$,
$1\leq k_0\leq n$, such that $|a_ia_{i+1}\cdots a_{i+(k_0-1)}|
\neq |a_ja_{j+1}\cdots a_{j+(k_0-1)}|$ $(a_{n+t}\equiv a_t)$.
Indeed, if there exist $i_0$ and $j_0$, $1\leq i_0<j_0\leq n$, such that
$|a_{i_0}\cdots a_{i_0+k}| = |a_{j_0}\cdots a_{j_0+k}|$ for all
$0\leq k\leq n-1$, then $|a_{i_0+k}| = |a_{j_0+k}|$ for all $0\leq
k\leq n-1$, that is, $|a_t| = |a_{t+(j_0-i_0)}|$ for all $1\leq
t\leq n$. This implies that $A$ has periodic weights, a
contradiction. Therefore, for any $1\leq i<j\leq n$, we have
$|a_ia_{i+1}\cdots a_{i+(k-1)}| \neq |a_ja_{j+1}\cdots
a_{j+(k-1)}|$ for some $k$. Note that $A^kA^{k*} =$ diag
$(\alpha_1, \ldots , \alpha_n)$ where $\alpha_t =
|a_ta_{t+1}\cdots a_{t+(k-1)}|^2$ for all $1\leq t\leq n$, and the
$(i,j)$-entry of $(A^kA^{k*})P$ and $P(A^kA^{k*})$ are
$\alpha_ip_{i j}$ and $p_{i j}\alpha_j$, respectively. Since
$(A^kA^{k*})P = P(A^kA^{k*})$ and all $a_j$'s are nonzero, hence
$\alpha_ip_{i j} = p_{i j}\alpha_j$ or $p_{i j} = 0$. From this
and $P = P^*$, we conclude that $P$ is a diagonal matrix, that is, $P
= $ diag $(p_{1 1}, \ldots , p_{n n})$ and $p_{j j} = 0$ or $1$
for all $1\leq j\leq n$. Moreover, $AP-PA = 0$ implies that
$a_1(p_{1 1}-p_{2 2})=a_2(p_{2 2}-p_{3 3})=\cdots
=a_n(p_{n n}-p_{1 1})=0$, since all $a_j$'s are nonzero, we deduce
that $P=0$ or $P=I_n$. This contradicts the fact that $A$ is
reducible. Therefore, $A$ has periodic weights.

\vspace{.5cm}

(b)$\Rightarrow$(c). Let $\tilde{A}$ be the $n$-by-$n$ weighted
shift matrix with weights $|a_1|, |a_2|, \ldots , |a_n|$. By
Lemma 2.1 (b), we obtain that $A$ is unitarily equivalent
to $e^{i\theta}\tilde{A}$, where $\theta =
\left(\sum_{i=1}^{n}\arg a_j\right)/n$. We want to construct an
$n$-by-$n$ unitary matrix $U$ such that $U^*\tilde{A}U = B\oplus
\omega B\oplus \cdots \oplus \omega^{m-1}B$ where $\omega =e^{2\pi
i/n}$, $m=n/k$ and $B$ is the $k$-by-$k$ weighted shift matrix
with weights $|a_1|, \ldots , |a_k|$. Let $V_j =$ diag $(1,
\omega^j, \omega^{2j}, \ldots , \omega^{(k-1)j})$ for $j = 1, 2,
\ldots , m-1$ and $U$ be the $n$-by-$n$ matrix
$$\frac{1}{\sqrt{m}}\left[\begin{array}{cccccccc}I_k & V_1 & V_2
& \cdots & V_{m-1} \\ I_k & \omega^kV_1 & (\omega^2)^kV_2 & \cdots
& (\omega^{m-1})^kV_{m-1} \\ \vdots & \vdots & \vdots & \ddots &
\vdots \\ I_k & \omega^{(m-1)k}V_1 & (\omega^2)^{(m-1)k}V_2 &
\cdots & (\omega^{m-1})^{(m-1)k}V_{m-1}\end{array}\right].$$ We
now check that $U$ is unitary. Indeed,  for any $1\leq s,t\leq m$,
the $(s,t)$-block of $UU^*$ is
\begin{eqnarray*}
&&\frac{1}{m}\left[\begin{array}{cccccc}I_k & \omega^{(s-1)k}V_1 & \omega^{2(s-1)k}V_2 & \cdots & \omega^{(m-1)(s-1)k}V_{m-1} \end{array}\right]\left[\begin{array}{cccccc}I_k \\ \bar{\omega}^{(t-1)k}V_1^* \\ \bar{\omega}^{2(t-1)k}V_2^* \\ \vdots \\ \bar{\omega}^{(m-1)(t-1)k}V_{m-1}^*\end{array}\right]\\
& =& \frac{1}{m}\left(\sum_{l=0}^{m-1}\omega^{l(s-t)k}\right)I_k.
\end{eqnarray*}
If $s=t$, then $\sum_{l=0}^{m-1}\omega^{l(s-t)k} = m$. On the
other hand, if $s \neq t$, then
$$\sum_{l=0}^{m-1}(\omega^{(s-t)k})^l =
\frac{1-\omega^{(s-t)km}}{1-\omega^{(s-t)k}} =
\frac{1-e^{2(s-t)\pi i}}{1-\omega^{(s-t)k}} = 0.$$ Thus
$UU^* = I_n$ or $U$ is unitary.

\vspace{.5cm}

Next, we check that $U(B\oplus \omega B\oplus \cdots \oplus
\omega^{m-1}B)U^* = \tilde{A}$. Write $B = B_1+B_2$ where $B_1$
(resp., $B_2$) is the $k$-by-$k$ weighted shift matrix with
weights $|a_1|, \ldots , |a_{k-1}|, 0$ (resp., $0, \ldots , 0,
|a_k|$). A simple computation shows that $V_jBV_j^* =
\omega^{-j}B_1+\omega^{(k-1)j}B_2$ for all $1\leq j\leq m-1$. For
any $1\leq s,t\leq m$, the $(s,t)$-block of $U(B\oplus \cdots \oplus
\omega^{m-1}B)U^*$ is
\begin{eqnarray*}
&& \frac{1}{m}\left[B \ \ \omega^{(s-1)k}V_1(\omega B) \
 \ \cdots \ \
\omega^{(m-1)(s-1)k}V_{m-1}(\omega^{m-1}B)\right]\left[\begin{array}{c}I_k \\
\bar{\omega}^{(t-1)k}V_1^* \\
\vdots \\ \bar{\omega}^{(m-1)(t-1)k}V_{m-1}^*\end{array}\right] \\
&=&\frac{1}{m}\left(B+\omega^{(s-t)k+1}V_1BV_1^*+\omega^{2((s-t)k+1)}V_2BV_2^*+\cdots
+\omega^{(m-1)((s-t)k+1)}V_{m-1}BV_{m-1}^*\right) \\
&=&\frac{1}{m}\left(\sum_{l=0}^{m-1}\omega^{l(s-t)k}\right)B_1+\frac{1}{m}\left(\sum_{l=0}^{m-1}\omega^{l(s-t+1)k}\right)B_2\\
&=&\frac{1}{m}\left(\sum_{l=0}^{m-1}e^{2l(s-t)\pi
i/m}\right)B_1+\frac{1}{m}\left(\sum_{l=0}^{m-1}e^{2l(s-t+1)\pi
i/m}\right)B_2 \\
&=&\left\{\begin{array}{ll}B_1 \hspace{5mm}& \mbox{if }  1\leq s= t\leq
m, \\ B_2 & \mbox{if }  (s,t) =(m,1) \mbox{ or } t=s+1, 1\leq s\leq m-1,  \\
0 & \mbox{otherwise.}  \end{array}\right.
\end{eqnarray*}
>From above, we obtain $$U(B\oplus \omega B\oplus \cdots \oplus
\omega^{m-1}B)U^* = \left[\begin{array}{cccc}B_1 & B_2 & & \\ &
B_1 & \ddots & \\ & & \ddots & B_{2} \\ B_2 & & &
B_1\end{array}\right] = \tilde{A}$$ as required.

\vspace{.5cm}

(c)$\Rightarrow$(a). This implication is trivial.
\hfill$\blacksquare$

\vspace{.5cm}

Let $A$ be an $n$-by-$n$ weighted
shift matrix with weights $a_1,\ldots , a_{n}$. The next proposition shows that
if $A$ has exactly one zero weight, then $A$ is irreducible.

\vspace{.5cm}

{\bf Proposition 3.2.} \emph{Let $A$ be an $n$-by-$n$ weighted
shift matrix with weights $a_1,\ldots , a_{n-1}, 0$, where
$a_j\neq 0$ for all $1\leq j\leq n-1$. Then $A$ is irreducible.}

\vspace{.5cm}

{\em Proof.} Let $\{e_1, \ldots , e_n\}$ be the standard basis for
$\mathbb{C}^n$, and $M$ be a nontrivial reducing subspace of $A$.
We want to show that $M = \mathbb{C}^n$. Indeed, let $x = [x_1\
\ldots \  x_n]^T$ be a nonzero vector in $M$, and $x_{j_0}$ be the
first nonzero entry of $x$. Then ${A^*}^{(n-j_0)}x
=\bar{a}_{j_0}\bar{a}_{j_0+1}\cdots \bar{a}_{n-1}x_{j_0}e_n$, it
follows that $e_n\in M$. Consequently, we have $A^je_n\in M$ for
all $j = 1, \ldots , n-1$. Since $A^je_n=a_{n-j}a_{n-j+1}\cdots a_{n-1}e_{n-j}$ for $j=1, \ldots, n-1$, hence $\{e_1, \ldots , e_n\}\subseteq M$ and $A$ is irreducible. \hfill$\blacksquare$

\vspace{.5cm}

We now give a complete characterization of $n$-by-$n$ weighted
shift matrices $A$ which are reducible.

\vspace{.5cm}

{\bf Corollary 3.3.} \emph{Let $A$ be an $n$-by-$n$} $(n\geq 2)$
\emph{weighted shift matrix with weights $a_1,\ldots , a_n$. Then
$A$ is reducible if and only if one of the following cases hold:}

(a) \emph{$a_i = a_j = 0$ for some $1\leq i<j \leq n$,}

(b) \emph{$A$ has periodic weights.}

\vspace{.5cm}

{\em Proof}. Assume that $A$ is reducible. From Proposition 3.2, we infer that either there are at least two weights of $A$ being zero, or all $a_j$'s are nonzero. Therefore, $A$ is either in case (a) or in case (b) from Theorem 3.1.

\vspace{.5cm}

To prove the converse, we first assume that $a_i=a_j=0$ for some $i, j$, $1\le i<j\le n$. By Lemma 2.1 (a), we may assume that $j=n$ and $1<i<n$. Then $A=A_1\oplus A_2$, where $A_1$ and $A_2$ are the weighted shift matrices with weights $a_1, \ldots, a_{i-1}, 0$ and $a_{i+1}, \ldots, a_{n-1}, 0$, respectively. This shows that $A$ is reducible. For case (b), if all $a_j$'s are nonzero, then our assertion follows from Theorem 3.1.  If $a_j=0$ for some $j$, since $A$ has periodic weights, then $a_{k+j}=a_j=0$ for some $k$, $1\le k\le\lfloor n/2\rfloor$. Hence $A$ is reducible from case (a). This completes the proof. \hfill$\blacksquare$

\vspace{.5cm}

We conclude this section by remarking that \cite[Theorem 1 (a)]{4} is an immediate consequence of Theorem 3.1.

\vspace{1cm}

\centerline {\bf\large 4. Numerical ranges}

\vspace{.4cm}

Recall that the \emph{numerical range} of an $n$-by-$n$ matrix $A$
is by definition the set $W(A) = \{\langle Ax,x\rangle : x\in
\mathbb{C}^n, \|x\| = 1 \}$, where $\langle \cdot , \cdot\rangle$
and $\|\cdot \|$ denote, respectively, the standard inner product
and Euclidean norm in $\mathbb{C}^n$. The numerical range is a
nonempty compact convex subset of the complex plane. It is
invariant under unitary equivalence and contains the eigenvalues.
For other properties of the numerical range, the reader can
consult \cite[Chapter 1]{3}.

\vspace{.4cm}

In recent years, properties of the numerical ranges of weighted
shift matrices have been intensely studied (cf. \cite{4, 5, 6}).
It was obtained that the numerical range of an $n$-by-$n$ weighted
shift matrix $A$ has the $n$-symmetry property, that is, $W(A) =
e^{2\pi i/n}W(A)$ (cf. \cite[Theorem 2.3]{1}). Moreover, if
$A$ has at least one zero weight, then $W(A)$ is a circular disc
centered at the origin. In this case, Stout \cite{6} gave a formula for the radius of the circular disc $W(A)$.
His formula involves the circularly
symmetric functions. In this section, we
will give the expansion of the Kippenhahn polynomial of $A$ in terms of the
circularly symmetric functions. Therefore, here we give a brief review of
the circularly symmetric functions, following Stout \cite{6}.

\vspace{.4cm}

Let $a_1, \ldots , a_n$ be complex numbers and $r$ be a
nonnegative integer. $S_0$ is defined to be $1$, while  for $r\geq
1$, $S_r(a_1, \ldots , a_n) = \sum\left\{ \prod^{r}_{k=1}a_{\pi
(k)} | \pi : (1,\ldots ,r)\rightarrow (1,\ldots ,n)\right.$, where
$\pi (k)+1<\pi (k+1)$ for $1\leq k < r$, and if $\pi (1)=1$ then
$\left.\pi (r)\neq n\right\}$. These have a nice description:
imagine a regular $n$-gon with vertices labeled $a_1$ through
$a_n$. Draw a convex $r$-gon in it, with vertices among the $a_j$
with the restriction that it can not use an edge of the original
polygon. Each term in $S_r(a_1, \ldots ,a_n)$ is the product of
the vertices of such an $r$-gon.

\vspace{.4cm}

These functions satisfy many identities. By \cite[P. 496]{6}, we
have the following: \\
(4.1) $S_1(a_1, \ldots , a_n) = \sum_{k=1}^{n}a_k$ if $n>1$,\\
(4.2) $S_{r}(a_{1},\ldots , a_{n})=0$ if $r>n/2$, \\
(4.3) $S_{r}(a_{1},\ldots , a_{n}) = S_{r}(a_2, \ldots , a_{n}, a_1 )$, \\
(4.4) $S_{r}(a_{1},\ldots , a_{n}, 0) = S_{r}(a_{1},\ldots ,a_{n}, 0, 0) $, \\
(4.5) $S_{r+1}(a_{1},\ldots , a_{n+1}, 0)= S_{r+1}(a_{1},\ldots ,
a_{n}, 0) + a_{n+1}S_{r}(a_{1},\ldots , a_{n-1},0)$.

\vspace{.4cm}

In particular, we need the following identities.

\vspace{.4cm}

{\bf Proposition 4.1.} \emph{Let $S_r(a_1, \ldots , a_n)$ be the
circularly symmetric function defined as above. Then}

(a) $S_{r}(a_{1},\ldots , a_{n}, 0)= S_{r}(a_{2},\ldots , a_{n},
0) + a_{1}S_{r-1}(a_{3},\ldots , a_{n}, 0)$, \emph{and}

(b) $S_{r}(a_{1},\ldots , a_{n})= S_{r}(a_{1},\ldots , a_{n}, 0) -
a_{1}a_{n}S_{r-2}(a_{3},\ldots , a_{n-2}, 0)$.

\vspace{.4cm}

{\em Proof.} By the definition of the circularly symmetric
function, it is clear that
\begin{eqnarray*}
&&S_r(a_1, \ldots , a_n, 0) \\
&=& \sum\left\{ \prod^{r}_{k=1}a_{\pi (k)} |  \pi : (1,\ldots
,r)\rightarrow (1,\ldots ,n), \mbox{  where  } \pi (k)+1<\pi (k+1)
\mbox{  for  } 1\leq k < r\right\}.
\end{eqnarray*}
(a) Note that
\begin{eqnarray*}
&&\left\{\pi : (1,\ldots ,r)\rightarrow (1,\ldots ,n) | \pi
(k)+1<\pi (k+1) \mbox{  for  } 1\leq k < r \right\} \\
&=&\left\{\pi : (1,\ldots ,r)\rightarrow (1,\ldots ,n) | \pi
(k)+1<\pi (k+1) \mbox{  for  } 1\leq k < r \mbox{  and  } \pi (1) = 1\right\} \\
&& \cup \left\{\pi : (1,\ldots ,r)\rightarrow (1,\ldots ,n) | \pi
(k)+1<\pi (k+1) \mbox{  for  } 1\leq k < r \mbox{  and  } \pi (1) \neq 1\right\} \\
&=&\left\{\pi : (1,\ldots ,r-1)\rightarrow (3,\ldots ,n) | \pi
(k)+1<\pi (k+1) \mbox{  for  } 1\leq k < r-1\right\} \\
&& \cup \left\{\pi : (1,\ldots ,r)\rightarrow (2,\ldots ,n) | \pi
(k)+1<\pi (k+1) \mbox{  for  } 1\leq k < r\right\}.
\end{eqnarray*}
Hence $S_r(a_1,\ldots ,a_n, 0) = S_{r}(a_{2},\ldots , a_{n}, 0) +
a_{1}S_{r-1}(a_{3},\ldots , a_{n}, 0)$.

\vspace{.4cm}

(b) By definition, we have
\begin{eqnarray*}
&&S_r(a_1,\ldots , a_n, 0)-S_r(a_1,\ldots ,a_n)\\
&=&\sum\left\{ \prod^{r}_{k=1}a_{\pi (k)} |  \pi : (1,\ldots
,r)\rightarrow (1,\ldots ,n), \pi (k)+1<\pi (k+1), 1\leq k < r,
\pi (1) = 1, \pi (r) =
n\right\} \\
&=&\sum\left\{ a_1a_n\prod^{r-1}_{k=2}a_{\pi (k)} |  \pi :
(2,\ldots ,r-1)\rightarrow (3,\ldots ,n-2), \pi
(k)+1<\pi (k+1), 2\leq k < r-1\right\} \\
&=&a_1a_n\sum\left\{\prod^{r-2}_{k = 1}a_{\alpha (k)} |
\alpha : (1,\ldots ,r-2)\rightarrow (3,\ldots ,n-2), \alpha (k)+1<\alpha (k+1), 1\leq k < r-2\right\} \\
&=&a_1a_nS_{r-2}(a_{3},\ldots , a_{n-2}, 0)
\end{eqnarray*}
as asserted. \hfill$\blacksquare$

\vspace{.4cm}

For an $n$-by-$n$ matrix $A$, consider the degree-$n$ homogeneous
polynomial $p_A(x,y,z) = \det(x\re A+y\im A+zI_n)$, where $\re A =
(A+A^*)/2$ and $\im A = (A-A^*)/(2i)$. A result of Kippenhahn
\cite[Theorem 10]{7} says that the numerical range $W(A)$ equals the
convex hull of the real points in the dual of the curve
$p_A(x,y,z) = 0$, that is, $$W(A) = \left\{a+ib\in\mathbb{C} : a,b
\mbox{  real,  }ax+by+z = 0 \mbox{  tangent to }p_A(x,y,z) =
0\right\}^\wedge .$$ Here, for any subset $\bigtriangleup$ of
$\mathbb{C}$, $\bigtriangleup^\wedge$ denote its convex hull, that
is, $\bigtriangleup^\wedge$ is the smallest convex set containing
$\bigtriangleup$. Therefore, the numerical range $W(A)$ is
completely determined by the Kippenhahn polynomial $p_A(x,y,z)$.

\vspace{.4cm}

The next theorem gives the expansion of the Kippenhahn polynomial
$p_A(x,y,z)$ of an $n$-by-$n$ weighted shift matrix $A$ in terms
of its weights.

\vspace{.4cm}

{\bf Theorem 4.2.} {\em Let $A$ be an $n$-by-$n$ $(n\geq3)$
weighted shift matrix with weights $a_1, \ldots, a_n$}. Then
\begin{eqnarray*}
p_{A}(x,y,z)
&=& z^n + \sum^{\lfloor \frac{n}{2}\rfloor}_{r=1} S_{r}(|a_{1}|^{2},\ldots ,|a_{n}|^{2})(x^2+y^2)^{r}(-\frac{1}{4})^{r}z^{n-2r} \\
&& + \frac{(-1)^{n+1}}{2^{n}}\left((x-iy)^n\prod^{n}_{j=1}a_{j} +
(x+iy)^n\prod^{n}_{j=1}\bar{a}_{j} \right).
\end{eqnarray*}

\vspace{.4cm}

For the proof of Theorem 4.2, we need the next two lemmas. The first lemma
is an immediate consequence of the result of Stout \cite[Lemma
1]{6}.

\vspace{.4cm}

{\bf Lemma 4.3.} \emph{Let $A$ be an $n$-by-$n$ weighted shift
matrix with weights $a_1, \ldots , a_{n-1}, 0$. Then }$$\det
(zI_n+\re A) = \sum_{r=0}^{\lfloor n/2 \rfloor}S_r(|a_1|^2, \ldots
, |a_{n-1}|^2, 0)(-\frac{1}{4})^rz^{n-2r}.$$

\vspace{.4cm}

{\bf Lemma 4.4.} \emph{Let $A$ be an $n$-by-$n$ weighted shift
matrix with weights $a_1, \ldots , a_{n-1}, 0$. Then } $$\det
(x\re A+y\im A+zI_n) = \sum_{r=0}^{\lfloor n/2
\rfloor}S_r(|a_1|^2, \ldots , |a_{n-1}|^2,
0)(x^2+y^2)^r(-\frac{1}{4})^rz^{n-2r}.$$

\vspace{.4cm}

{\em Proof.} Let $\tilde{A}$ be the $n$-by-$n$ weighted shift
matrix with weights $a_1(x-iy), a_2(x-iy), \ldots , a_{n-1}(x-iy),
0$. It is easily seen that $\re \tilde{A} = x\re A+y\im A$.
Therefore, by Lemma 4.3, we have
\begin{eqnarray*}
&&\det (x\re A+y\im A+zI_n) \\&=&\det (zI_n+\re \tilde{A}) \\
&=&\sum_{r=0}^{\lfloor n/2 \rfloor}S_r(|a_1|^2(x^2+y^2), \ldots ,
|a_{n-1}|^2(x^2+y^2), 0)(-\frac{1}{4})^rz^{n-2r}\\
&=&\sum_{r=0}^{\lfloor n/2 \rfloor}S_r(|a_1|^2, \ldots ,
|a_{n-1}|^2, 0)(x^2+y^2)^r(-\frac{1}{4})^rz^{n-2r}
\end{eqnarray*}
as asserted. \hfill$\blacksquare$

\vspace{.4cm}

We are now ready to prove Theorem 4.2. For simplicity, let
$A[i_1,\ldots , i_m]$ denote the $(n-m)$-by-$(n-m)$ principal
submatrix of $A$ obtained by deleting its rows and columns indexed
by $i_1,\ldots , i_m$.

\vspace{.4cm}

{\em Proof of Theorem} 4.2. We now expand the determinant of
$x\re A + y\im A + zI_{n}$ by minor along its $n$th column to
obtain
\begin{eqnarray*}
&&\det (x\re A + y\im A + zI_{n}) \\
&=& z\det (C_1+zI_{n-1})+\frac{a_{n-1}}{2}(x-iy)(-1)^{2n-1}
d_{n-1,n} + (-1)^{n+1}\frac{\bar{a}_{n}}{2}(x+iy)d_{1,n},
\end{eqnarray*}
where $C_1=x\re A[n]+y\im A[n]$ and $(-1)^{n+j}d_{j,n}$ denotes the
cofactor of the $(j,n)$-entry of $x\re A + y\im A + zI_{n}$ in
$x\re A + y\im A + zI_{n}$, $j = 1,n-1$. The expansion of the
determinant $d_{1,n}$ (resp., $d_{n-1,n}$) along its last
row (resp., its last column) yields $$d_{1,n} = \bar{a}_1\cdots
\bar{a}_{n-1}(\frac{x+iy}{2})^{n-1}+(-1)^n\frac{a_n}{2}(x-iy)\det
(C_2+zI_{n-2})$$
$$(\mbox{resp}., d_{n-1,n} = (-1)^na_1\cdots
a_{n-2}a_n(\frac{x-iy}{2})^{n-1}+\frac{\bar{a}_{n-1}}{2}(x+iy)\det
(C_3+zI_{n-2})),$$ where $C_2 = x\re A[1,n]+y\im A[1,n]$ (resp.,
$C_3=x\re A[n-1,n]+y\im A[n-1,n]$).
Hence
\begin{eqnarray*}
&&\det (x\re A+y\im A+zI_n) \\
&=&z\det (C_1+zI_{n-1}) -\frac{|a_{n-1}|^2}{4}(x^2+y^2)\det
(C_3+zI_{n-2})-\frac{|a_n|^2}{4}(x^2+y^2)\det
(C_2+zI_{n-2})\\
&& + (-1)^{n+1}\bar{a}_1\cdots
\bar{a}_n(\frac{x+iy}{2})^n+(-1)^{n+1}a_1\cdots
a_{n}(\frac{x-iy}{2})^n.
\end{eqnarray*}
By Lemma 4.4, we obtain that
$$z\det (C_1+zI_{n-1}) = \sum_{r=0}^{\lfloor (n-1)/2\rfloor}S_r(|a_1|^2,\ldots , |a_{n-2}|^2, 0)(x^2+y^2)^r(\frac{-1}{4})^rz^{n-2r}.$$
Moreover,
\begin{eqnarray*}
&&-\frac{|a_n|^2}{4}(x^2+y^2)\det (C_2+zI_{n-2})\\& =
&\sum_{r=0}^{\lfloor (n-2)/2\rfloor}|a_n|^2S_r(|a_2|^2,\ldots ,
|a_{n-2}|^2, 0)(x^2+y^2)^{r+1}(\frac{-1}{4})^{r+1}z^{n-2-2r}\\
&=&\sum_{r=1}^{\lfloor n/2\rfloor}|a_n|^2S_{r-1}(|a_2|^2,\ldots ,
|a_{n-2}|^2, 0)(x^2+y^2)^{r}(\frac{-1}{4})^{r}z^{n-2r},
\end{eqnarray*}
and
\begin{eqnarray*}
&&-\frac{|a_{n-1}|^2}{4}(x^2+y^2)\det (C_3+zI_{n-2})\\
&=&\sum_{r=0}^{\lfloor
(n-2)/2\rfloor}|a_{n-1}|^2S_r(|a_1|^2,\ldots ,
|a_{n-3}|^2,0)(x^2+y^2)^{r+1}(\frac{-1}{4})^{r+1}z^{n-2-2r}\\
&=&\sum_{r=1}^{\lfloor
n/2\rfloor}|a_{n-1}|^2S_{r-1}(|a_1|^2,\ldots , |a_{n-3}|^2,
0)(x^2+y^2)^{r}(\frac{-1}{4})^{r}z^{n-2r}.
\end{eqnarray*}
Note that if $n$ is odd, then $\lfloor (n-1)/2\rfloor = \lfloor
n/2\rfloor$. On the other hand, if $n$ is even, then $n/2>(n-1)/2$ and, by
Equation (4.2), we have $$0 = S_{n/2}(|a_1|^2,\ldots , |a_{n-2}|^2, 0) =
S_{(n/2)-1}(|a_2|^2,\ldots , |a_{n-2}|^2, 0) =
S_{(n/2)-1}(|a_1|^2,\ldots , |a_{n-3}|^2, 0).
$$ Therefore, we deduce that
\begin{eqnarray*}
&&\det (x\re A+y\im A+zI_n) \\
&=& z^n + \sum^{\lfloor
n/2\rfloor}_{r=1}
\alpha_r(x^2+y^2)^{r}(-\frac{1}{4})^{r}z^{n-2r}
+\frac{(-1)^{n+1}}{2^{n}}\left((x-iy)^n\prod^{n}_{j=1}a_{j} +
(x+iy)^n\prod^{n}_{j=1}\bar{a}_{j} \right),
\end{eqnarray*}
where
\begin{eqnarray*}
\alpha_r &=& S_{r}(|a_1|^2,\ldots , |a_{n-2}|^2,
0)+|a_{n-1}|^2S_{r-1}(|a_1|^2,\ldots , |a_{n-3}|^2,
0)\\
&& +|a_n|^2S_{r-1}(|a_2|^2,\ldots , |a_{n-2}|^2, 0).
\end{eqnarray*}
for $r=1, \ldots, \lfloor n/2\rfloor$. We then apply
Equation (4.5) and Proposition 4.1 (a) and (b) to obtain
\begin{eqnarray*}
\alpha_r &=& S_{r}(|a_1|^2,\ldots , |a_{n-1}|^2,
0)+|a_n|^2S_{r-1}(|a_2|^2,\ldots , |a_{n-2}|^2, 0)\\
&=&S_{r}(|a_1|^2,\ldots , |a_{n-1}|^2,
0)+|a_n|^2S_{r-1}(|a_1|^2,\ldots , |a_{n-2}|^2,
0)\\
&&-|a_1|^2|a_n|^2S_{r-2}(|a_3|^2,\ldots , |a_{n-2}|^2, 0)\\
&=&S_{r}(|a_1|^2,\ldots , |a_{n}|^2,
0)-|a_1|^2|a_n|^2S_{r-2}(|a_3|^2,\ldots , |a_{n-2}|^2, 0)\\
&=&S_{r}(|a_1|^2,\ldots , |a_{n}|^2)
\end{eqnarray*}
for $r=1, 2, \ldots ,\lfloor n/2\rfloor$. This completes the proof.
\hfill$\blacksquare$

\vspace{.4cm}

We now restrict our attention to the Kippenhahn polynomial of a
weighted shift matrix. Let us recall some other known properties
of curves  in the complex projective plane $\mathbb{C}P^2$. Let
$p(x,y,z)$ be a degree-$n$ homogeneous polynomial and $\Gamma$ be
the dual curve of $p(x,y,z) = 0$. It is clear that $p(\alpha x,
\alpha y, \alpha z) = \alpha^np(x,y,z)$ for some scalar $\alpha$.
A point $\lambda = a+ib$, a, b real, is called a \emph{real focus}
of $\Gamma$ if $p(1, \pm i, -(a\pm ib)) = 0$ is satisfied.
Consequently, the eigenvalues of an $n$-by-$n$ matrix $T$ are
exactly the real foci of the dual curve of $p_T = 0$ (cf.
\cite[Theorem 11]{7}).
Moreover, for any $\theta\in \mathbb{R}$, since $\re
(e^{i\theta}T) = (\cos\theta) \re T - (\sin\theta) \im T$ and  $\im
(e^{i\theta}T) = (\cos\theta) \im T + (\sin\theta) \re T$, then
\begin{eqnarray*}
p_{e^{i\theta}T} (x,y,z)&=& \det \left((x\cos\theta+y\sin\theta)\re
T+(-x\sin\theta+y\cos\theta)\im T+zI_n\right)\\
&=&p_T(x\cos\theta+y\sin\theta, -x\sin\theta+y\cos\theta, z).
\end{eqnarray*}
Among other things, if $A$ and $B$ are $n$-by-$n$ matrices so that
$A$ is unitary equivalent to $B$, then $p_A(x,y,z) = p_B(x,y,z)$.

\vspace{.4cm}

Let $A$ be an $n$-by-$n$ weighted shift matrix with weights $a_1,
\ldots , a_n$. Lemma 2.1 (b) says that $A$ is unitarily
equivalent to $\omega_n^jA$ for all $1\leq j\leq n$, where
$\omega_n = e^{2\pi i/n}$. Thus $p_A(x,y,z) =
p_{\omega_n^jA}(x,y,z)$ for all $1\leq j\leq n$. Furthermore, if $A$
is reducible, Theorem 3.1 says that $A$ is unitarily equivalent to
$B\oplus (\omega_nB)\oplus \cdots \oplus (\omega_n^{(n/k)-1}B)$,
where $B$ is the $k$-by-$k$ weighted shift matrix with weights
$|a_1|e^{i\theta}, \ldots , |a_k|e^{i\theta}$, $k$ is a factor of
$n$ and $\theta = (\sum_{i=1}^{n}\arg a_j)/n$. Consequently, we
have
\begin{eqnarray*}
p_A(x,y,z)&=&\prod_{j=0}^{(n/k)-1}p_{\omega_n^jB}(x,y,z)\\
&=&\prod_{j=1}^{n/k}p_B(x\cos\theta_j+y\sin\theta_j,
-x\sin\theta_j+y\cos\theta_j, z)
\end{eqnarray*}
where $\theta_j = 2(j-1)\pi /n$ for $j = 1,\ldots , n/k$. For an
$n$-by-$n$ matrix $T$, we known that if $T$ is reducible then
$p_T$ is also reducible. But the converse is not true in general.
The next proposition shows that if $p_A$ is reducible then $p_A$
must be of the form described above.

\vspace{.4cm}

{\bf Proposition 4.5.} {\em Let $A$ be an $n$-by-$n$ $(n\geq 2)$
weighted shift matrix. If $p_A$ is reducible and has an
irreducible factor $q(x,y,z)$ of degree-$k$, then $n$ is divisible
by $k$, and}
$$p_A(x,y,z) = \prod_{j=1}^{n/k}q(x\cos\theta_j+y\sin\theta_j, -x\sin\theta_j+y\cos\theta_j, z),$$ where $\theta_j = 2(j-1)\pi /n$
for $j = 1,\ldots , n/k$.

\vspace{.4cm}

{\em Proof.} For abbreviation, we let $$q_j(x,y,z) =
q(x\cos\theta_j+y\sin\theta_j, -x\sin\theta_j+y\cos\theta_j, z)$$ for $j =
1,\ldots , n$. Note that $q_1(x,y,z) = q(x,y,z)$. We first check
that $q_j(x,y,z)$ is a factor of $p_A(x,y,z)$ for all $1\leq j\leq
n$. Indeed, since $q$ is a factor of $p_A$, then $p_{A}(x,y,z) =
q(x,y,z)r(x,y,z)$ for some homogeneous polynomial $r(x,y,z)$. Note
that $A$ is unitarily equivalent to $e^{i\theta_j}A$ for all
$1\leq j\leq n$, hence
\begin{eqnarray*}
&&p_A(x,y,z)\\
&=&p_{e^{i\theta_j}A}(x,y,z)\\
&=&p_{A}(x\cos\theta_j+y\sin\theta_j, -x\sin\theta_j+y\cos\theta_j, z)\\
&=& q(x\cos\theta_j+y\sin\theta_j, -x\sin\theta_j+y\cos\theta_j, z)r(x\cos\theta_j+y\sin\theta_j, -x\sin\theta_j+y\cos\theta_j, z),\\
&=& q_{j}(x,y,z)r(x\cos\theta_j+y\sin\theta_j,
-x\sin\theta_j+y\cos\theta_j, z)
\end{eqnarray*}
for all $1\leq j\leq n$. Therefore, $q_j$ is a factor of $p_A$ for
all $1\leq j\leq n$.

\vspace{.4cm}

Next, we want to show that $n$ is divisible by $k$. For each $j=1,
2, \ldots , n$, let $\Gamma_j$ be the dual curve of $q_j(x,y,z) =
0$ and $E_j$ be the set of all real foci of $\Gamma_j$. We claim
that $E_j = e^{i\theta_j}E_1$ for all $1\leq j\leq n$. Indeed, if
$\lambda = a+ib\in E_1$, where $a$ and $b$ are real, then $
e^{i\theta_j}\lambda = (a\cos \theta_j-b\sin \theta_j)+i(a\sin
\theta_j+b\cos \theta_j)$, and
\begin{eqnarray*}
&&q_j(1, \pm i, -\left((a\cos \theta_j-b\sin \theta_j)\pm i(a\sin
\theta_j+b\cos \theta_j)\right)) \\
&=&q(\cos \theta_j\pm i\sin \theta_j, -\sin \theta_j\pm i\cos
\theta_j, -(\cos \theta_j \pm i\sin\theta_j)(a\pm ib))\\
&=&(\cos \theta_j \pm \sin \theta_j)^k\cdot q(1, \pm i, -(a\pm
ib)) = 0.
\end{eqnarray*}
Note that every real focus of $\Gamma_j$ is an eigenvalue of $A$,
and every eigenvalue of $A$ is simple, it follows that all real
foci of $\Gamma_j$ are distinct. Hence $E_j = e^{i\theta_j}E_1$
for all $1\leq j\leq n$ as claimed. Moreover, for any $1\leq i\neq
j\leq n$, if $q_i\neq q_j$, since $p_A$ is divisible by $q_iq_j$
and the eigenvalue of $A$ is simple, it implies that
$E_i\cap E_j = \emptyset$. Therefore,
we conclude that either $E_i=E_j$ or $E_i\cap E_j = \emptyset$  for any $1\leq i<j \leq n$.

\vspace{.4cm}

Now, if $E_1\cap E_j = \emptyset$ for all $2\leq j\leq n$, we will
prove $E_s\cap E_t = \emptyset$ for any $1\leq s<t\leq n$. On the contrary, if $E_s=E_t$ for some $s$, $t$, since
$E_s = e^{i\theta_s}E_1$ and $E_t = e^{i\theta_t}E_1$, it follows
that $e^{i\theta_s}E_1 = e^{i\theta_t}E_1$ or $E_1 =
e^{i(\theta_t-\theta_s)}E_1 = e^{i\theta_{t-s+1}}E_1 = E_{t-s+1}$,
a contraction. Hence these sets $E_j$'s are disjoint. Note that
$E_j\subseteq \sigma (A)$ and $E_j$ contains
exactly $k$ elements for all $1\leq j\leq n$. Thus
$\bigcup_{j=1}^{n}E_j\subseteq \sigma (A)$ and $k\cdot n\leq n$.
This clearly forces that $k=1$ and $\sigma (A) =
\bigcup_{j=1}^{n}E_j$, that is, $p_A = \Pi_{j=1}^{n}q_j$ as
desired.

\vspace{.4cm}

On the other hand, if $E_1 = E_j$ for some $j$, $2\leq j\leq n$. Let
$j_{0} = \min \{j  :  E_{j} = E_{1}, 2\leq j\leq n\}$ and
$m=j_0-1$. Then $E_{m+1}=E_1$, in consequence, $E_{m+j} =
E_j$ for all  $2\leq j\leq n$, because $E_{m+j} =
e^{i\theta_{m+j}}E_1 = e^{2(m+j-1)\pi i/n}E_1 = e^{2(j-1)\pi
i/n}\cdot e^{2m\pi i/n}E_1 = e^{i\theta_j}\cdot
e^{i\theta_{m+1}}E_1 = e^{i\theta_j}E_{m+1} = e^{i\theta_j}E_1 =
E_j$. From this, we have $\bigcup_{j=1}^{n}E_j =
\bigcup_{j=1}^{m}E_j\subseteq \sigma (A)$. But $\sigma (A)
=\left\{e^{i\theta_1}\lambda, e^{i\theta_2}\lambda, \ldots ,
e^{i\theta_n}\lambda\right\}$, where $\lambda = (a_1\cdots
a_n)^{1/n}$. It follows that $\sigma (A) \subseteq
\bigcup_{j=1}^{n}e^{i\theta_j}E_1 = \bigcup_{j=1}^{n}E_j =
\bigcup_{j=1}^{m}E_j\subseteq \sigma(A)$. Hence $\sigma (A) =
\bigcup_{j=1}^{m}E_j$, $n = km$ and $p_A = \prod_{j=1}^{m}q_j$.
This completes the proof. \hfill$\blacksquare$

\vspace{.4cm}

Let $A$ and $B$ be $n$-by-$n$ matrices. By Kippenhahn's result
\cite[Theorem 10]{7}, we known that if $p_A = p_B$ then $W(A) =
W(B)$. But the converse is not true in general. \cite[Example
2.1]{8} gives a counterexample. For the converse, if $W(A) =
W(B)$, then \cite[Proposition 2.3]{8} says that $p_A$ and $p_B$
have a common irreducible factor. Moreover, if $p_A$ is reducible, then $W(A)=W(B)$ if and only if $p_A = p_B$ (cf.
\cite[Corollary 2.4]{8}). Now, if $A$ and $B$ are $n$-by-$n$
weighted shift matrices, the next theorem shows that $W(A) = W(B)$
if and only if $p_A = p_B$, even if $p_A$ and $p_B$ are reducible.

\vspace{.4cm}

{\bf Theorem 4.6.} \emph{Let $A$ and $B$ be $n$-by-$n$ weighted
shift matrices with weights $a_1, \ldots , a_n$ and $b_1, \ldots ,
b_n$, respectively. Then the following statements are equivalent:}

(a) $W(A) = W(B)$,

(b) $p_{A} = p_{B}$,

(c) $S_{r}(|a_{1}|^2,\ldots ,|a_{n}|^2) = S_{r}(|b_{1}|^2,\ldots
,|b_{n}|^2)$ \emph{ for all $1\leq r\leq \lfloor n/2
\rfloor$ and $a_1\cdots a_n = b_1\cdots b_n$.}

\vspace{.4cm}

{\em Proof.} The equivalence of (b) and (c) follows from Theorem 4.2.
The implication (b)$\Rightarrow$(a) is a consequence of
Kippehhahn's result. To prove (a)$\Rightarrow$(b), assume that (a)
holds. If $p_A$ is irreducible, then $p_A = p_B$ follows from
\cite[Corollary 2.4]{8}. If $p_A$ is reducible, \cite[Proposition
2.3]{8} says that $p_A$ and $p_B$ have a common irreducible factor
$q$. Hence our assertion follows form Proposition 4.5. This
completes the proof. \hfill$\blacksquare$

\vspace{.4cm}

Let $A$ be an $n$-by-$n$ weighted shift matrix with positive
weights $a_1, \ldots ,a_n$. Since all $a_j$'s are positive, thus
$W(A) = W(A^*)$. Let $B$ be an $n$-by-$n$ weighted shift matrix
with positive weights $b_1, \ldots ,b_n$. If $b_j = a_{k+j}$ for
all $j$, for some fixed $k$, then $B$ is unitarily equivalent to $A$
and $W(B) = W(A)$. On the other hand, if $b_j = a_{(n+1)-j}$ for
all $j$, then $B$ is unitarily equivalent to $A^*$ and $W(B) =
W(A^*) = W(A)$. For the converse, if $W(B) = W(A)$, it is natural
to ask whether $B$ is unitarily equivalent to $A$ or $A^*$. The
following example shows that this is not the case. This
example can be easily constructed by Theorem 4.6.

\vspace{.4cm}

{\bf Example 4.7.} Let $$A=\left[\begin{array}{ccc}0 & 1 &
\\ & 0 & \frac{\sqrt{2}}{2} \\ \frac{\sqrt{30}}{4} & &
0\end{array}\right] \ \text{and} \ B=\left[\begin{array}{ccc}0 &
\sqrt{2} &  \\ & 0 & \frac{\sqrt{3}}{2} \\ \frac{\sqrt{10}}{4} & &
0\end{array}\right].$$ Form Theorem 4.6, we need to check that
$\det A = \det B$ and $S_1(1^2, (\sqrt{2}/2)^2,
(\sqrt{30}/4)^2)=S_1((\sqrt{2})^2, (\sqrt{3}/2)^2,
(\sqrt{10}/4)^2)$. By direct computation, we have $\det A
= \sqrt{60}/4 = \det B$ and $S_1(1^2, (\sqrt{2}/2)^2,
(\sqrt{30}/4)^2)=2+(1/2)+(30/16)
=27/8=2+(3/4)+(10/16)=S_1((\sqrt{2})^2,
(\sqrt{3}/2)^2, (\sqrt{10}/4)^2)$. Therefore,
Theorem 4.6 yields that $W(A) = W(B)$. But $B$ is neither
unitarily equivalent to $A$ nor unitarily equivalent to $A^*$ from
Theorem 2.3.

\vspace{.4cm}

Furthermore, let $A$ and $B$ be $n$-by-$n$ reducible weighted
shift matrices with positive weights. One may ask whether $B$ is
unitarily equivalent to $A$ or $A^*$ if $W(A) = W(B)$. The next
example shows that the answer is negative.

\vspace{.4cm}

{\bf Example 4.8.} Let $A$ and $B$ be the $3$-by-$3$ weighted
shift matrices as in Example 4.7, respectively, and $\tilde{A}$
and $\tilde{B}$ be $6$-by-$6$ weighted shift matrices with
weights$1, \sqrt{2}/2$, \\$\sqrt{30}/4, 1, \sqrt{2}/2,
\sqrt{30}/4$ and $\sqrt{2}, \sqrt{3}/2, \sqrt{10}/4, \sqrt{2},
\sqrt{3}/2, \sqrt{10}/4$, respectively. By Theorem 3.1,
$\tilde{A}$ (resp., $\tilde{B}$) is unitarily equivalent to
$A\oplus (e^{\pi i/3}A)$ (resp., $B\oplus (e^{\pi i/3}B)$).
Example 4.7 yields that $W(\tilde{A}) = W(\tilde{B})$. But $\tilde{B}$ is
neither unitarily equivalent to $\tilde{A}$ nor unitarily
equivalent to $\tilde{A^*}$ from Theorem 2.3.

\end{document}